# ON SURROGATE DIMENSION REDUCTION FOR MEASUREMENT ERROR REGRESSION: AN INVARIANCE LAW

BY BING LI[1] AND XIANGRONG YIN

*Pennsylvania State University and University of Georgia*

We consider a general nonlinear regression problem where the predictors contain measurement error. It has been recently discovered that several well-known dimension reduction methods, such as OLS, SIR and pHd, can be performed on the surrogate regression problem to produce consistent estimates for the original regression problem involving the unobserved true predictor. In this paper we establish a general invariance law between the surrogate and the original dimension reduction spaces, which implies that, at least at the population level, the two dimension reduction problems are in fact equivalent. Consequently we can apply *all* existing dimension reduction methods to measurement error regression problems. The equivalence holds exactly for multivariate normal predictors, and approximately for arbitrary predictors. We also characterize the rate of convergence for the surrogate dimension reduction estimators. Finally, we apply several dimension reduction methods to real and simulated data sets involving measurement error to compare their performances.

**1. Introduction.** We consider dimension reduction for regressions in which the predictor contains measurement error. Let $X$ be a $p$-dimension random vector representing the true predictor and $Y$ be a random variable representing the response. In many applications we cannot measure $X$ (e.g., blood pressure) accurately, but instead observe a surrogate $r$-dimensional predictor $W$ that is related to $X$ through the linear equation

$$W = \gamma + \Gamma^T X + \delta, \qquad (1)$$

where $\gamma$ is an $r$-dimensional nonrandom vector, $\Gamma$ is a $p$ by $r$ nonrandom matrix and $\delta$ is an $r$-dimensional random vector independent of $(X, Y)$. The

Received June 2005; revised November 2006.
[1]Supported in part by NSF Grants DMS-02-04662 and DMS-04-05681.
*AMS 2000 subject classifications.* 62G08, 62H12.
*Key words and phrases.* Central spaces, central mean space, invariance, regression graphics, surrogate predictors and response, weak convergence in probability.

goal of the regression analysis is to find the relation between the response $Y$ and the true, but unobserved, predictor $X$. This type of regression problem frequently occurs in practice and has been the subject of extensive studies, including, for example, those that deal with linear models (Fuller [15]), generalized linear models (Carroll [2], Carroll and Stefanski [5]), nonlinear models (Carroll, Ruppert and Stefanski [4]) and nonparametric models (see Pepe and Fleming [23]).

Typically there is an auxiliary sample which provides information about the relation between the original predictor $X$ and the surrogate predictor $W$, for example, by allowing us to estimate $\Sigma_{WX} = \text{cov}(W, X)$. Using this covariance estimate we can adjust the surrogate predictor $W$ to align it as much as possible with the true predictor $X$. At the population level this is realized by regressing $W$ on $X$, that is, adjusting $W$ to $U = \Sigma_{XW}\Sigma_W^{-1}W$, where $\Sigma_{XW} = \text{cov}(X, W)$ and $\Sigma_W = \text{var}(W)$. The fundamental question that will be answered in this paper is this: If we perform a dimension reduction operation on the surrogate regression problem of $Y$ versus $U$, will the result correctly reflect the relation between $Y$ and the true predictor $X$?

In the classical setting where the true predictor $X$ is observed, the dimension reduction problem can be briefly outlined as follows. Suppose that $Y$ depends on $X$ only through a lower dimensional vector of linear combinations of $X$, say $\beta^T X$, where $\beta$ is a $p$ by $d$ matrix with $d \leq p$. Or more precisely, suppose that $Y$ is independent of $X$ conditioning on $\beta^T X$, which will be denoted by

$$(2) \qquad Y \perp\!\!\!\perp X | \beta^T X.$$

The goal of dimension reduction is to estimate the directions of column vectors of $\beta$, or the column space of $\beta$. Note that the above relation will not be affected if $\beta$ is replaced by $\beta A$ for any nonsingular $p \times p$ matrix $A$. This is why the column space of $\beta$, rather than $\beta$ itself, is the object of interest in dimension reduction. A dimension reduction space provides us with a set of important predictors among all the linear combinations of $X$, with which we could perform exploratory data analysis or finer regression analysis without having to fit a nonparametric regression over a large number of predictors. Classical estimators of the dimension reduction space include ordinary least square (OLS) (Li and Duan [21], Duan and Li [13]), sliced inverse regression (SIR) (Li [19]), principle Hessian directions (pHd) (Li [20]) and the sliced inverse variance estimators (SAVE) (Cook and Weisberg [11]).

It has been discovered that some of these dimension reduction methods can be performed on the adjusted surrogate predictor $U$ to produce consistent estimates of at least some vectors in the column space of $\beta$ in (2) that describes the relation between $Y$ and the (unobserved) true predictor $X$. The first paper in this area is Carroll and Li [3], which demonstrated this phenomenon for OLS and SIR, and introduced the corresponding estimators of $\beta$ in the measurement error context. More recently, Lue [22] established that the pHd method, when applied to the surrogate problem $(U, Y)$, also yields consistent estimators of vectors in the column space of $\beta$. This work opens up the possibility of using available dimension reduction techniques to estimate $\beta$ by simply pretending $U$ is the true predictor $X$.

In this paper we will establish a general equivalence between the dimension reduction problem of $Y$ versus $U$ and that of $Y$ versus $X$. That is,

(3) $\qquad Y \perp\!\!\!\perp X | \beta^T X \quad \text{if and only if} \quad Y \perp\!\!\!\perp U | \beta^T U.$

This means that dimension reduction for the surrogate regression problem of $Y$ versus $U$ and that for the original regression problem of $Y$ versus $X$ are in fact equivalent at the population level. Thus the phenomena discovered by the above work are special cases of a very general invariance pattern—we can, in fact, apply *any* consistent dimension reduction method to the surrogate regression problem of $Y$ versus $U$ to produce consistent dimension reduction estimates for the original regression problem of $Y$ versus $X$. This fundamental relation is of practical importance, because OLS, SIR and pHd have some well-known limitations. For example, SIR does not perform well when the regression surface is symmetric about the origin, and pHd does not perform well when the regression surface lacks a clear quadratic pattern (or what is similar to it). New methods have recently been developed that can, in different respects and to varying degrees, remedy these shortcomings; see, for example, Cook and Li [9, 10], Xia et al. [25], Fung et al. [16], Yin and Cook [27, 28] and Li, Zha and Chiaromonte [18]. This equivalence allows us to choose among the broader class of dimension reduction methods to tackle the difficult situations in which the classical methods become inaccurate.

Sometimes the main purpose of the regression analysis is to infer the conditional mean $E(Y|X)$ or more generally conditional moments such as $E(Y^k|X)$. For example, in generalized linear models we are mainly interested in estimating the conditional mean $E(Y|X)$, and for regression with heteroscedasticity we may be interested in both the conditional mean $E(Y|X)$ and the conditional variance $\text{var}(Y|X)$. In these cases it is sensible to treat the conditional moments such as $E(Y|X)$ and $\text{var}(Y|X)$ as the objects of interest and the rest of the conditional distribution $f(Y|X)$ as the (infinite dimensional) nuisance parameter, and reformulate the dimension reduction problem to reflect this hierarchy. This was carried out in Cook and Li [9] and Yin and Cook [26], which introduced the notions of the central mean space and central moment space as well as methods to estimate them. If there is a $p$ by $d$ matrix $\beta$ with $d \leq p$ such that $E(Y|X) = E(Y|\beta^T X)$, then we call the column space of $\beta$ a dimension reduction space for the conditional mean $E(Y|X)$. More generally, the dimension reduction space for the $k$th conditional moment $E(Y^k|X)$ is defined as above with $Y$ replaced by $Y^k$. In this paper we will also establish the equivalence between the dimension reduction spaces for the $k$-conditional moments of the surrogate and the original regressions. That is,

(4) $\quad E(Y^k|X) = E(Y^k|\beta^T X) \quad \text{if and only if} \quad E(Y^k|U) = E(Y^k|\beta^T U).$

The above invariance relations will be shown to hold exactly under the assumption that $X$ and $\delta$ are multivariate normal; a similar assumption was also used in Carroll and Li [3] and Lue [22]. For arbitrary predictor and measurement error, we will establish an approximate invariant relation. This is based on the fact that, when $p$ is modestly large, most projections of a random vector are approximately normal (Diaconis and Freedman [12], Hall and Li [17]). Simulation studies indicate that the approximate invariance law holds for surprisingly small $p$ (as small as 6) and for severely nonnormal predictors.

This paper will be focused on the dimension reduction problems defined through relationships such as $Y \perp\!\!\!\perp X|\beta^T X$. A more general problem can be formulated as $Y \perp\!\!\!\perp X|t(X)$, where $t(X)$ is a (possibly nonlinear) function of $X$; see Cook [8]. Surrogate dimension reduction in this general sense is not covered by this paper, and remains an important open problem.

In Section 2 we introduce some basic issues and concepts related to measurement error problems and dimension reduction, as well as some machinery that will be repeatedly used in our further exposition. Equivalence (3) will be established in Section 3 for the case where $\Gamma$ in (1) is a $p$ by $p$ nonsingular matrix. Equivalence (3) for general $\Gamma$ will be shown in Section 4. In Section 5 we will establish equivalence (4). The approximate equivalence for general predictors and measurement errors will be developed in Section 6. In Section 7 we will turn our attention to a general estimation procedure for surrogate dimension reduction and study its convergence rate. In Section 8 we conduct a simulation study to compare different surrogate dimension reduction methods. In Section 9 we apply the invariance law to analyze a managerial behavior data set (Fuller [15]) that involves measurement errors. Some technical results will be proved in the Appendix.

**2. Preliminaries.** In this section we lay out some basic concepts and notation. For a pair of random vectors $V_1$ and $V_2$, we will use $\Sigma_{V_1 V_2}$ to denote the covariance matrix $\text{cov}(V_1, V_2)$, and for a random vector $V$, we will use $\Sigma_V$ to denote the variance matrix $\text{var}(V)$. If $V_1$ and $V_2$ are independent, then we write $V_1 \perp\!\!\!\perp V_2$; if they are independent conditioning on a third random element $V_3$, then we write $V_1 \perp\!\!\!\perp V_2|V_3$. If a matrix $\Sigma$ is positive definite, then we write $\Sigma > 0$. For a matrix $A$, the space spanned by its columns will be denoted by $\text{span}(A)$. If a matrix $A$ has columns $a_1, \ldots, a_p$, then $\text{vec}(A)$ denotes the vector $(a_1^T, \ldots, a_p^T)^T$. If $A, B, C$ are matrices, then $\text{vec}(ABC) = (C^T \otimes A)\text{vec}(B)$, where $\otimes$ denotes the tensor product between matrices.

In a measurement error problem we observe a *primary* sample on $(W, Y)$, which allows us to study the relation between $Y$ and $W$, and an *auxiliary* sample that allows us to estimate $\Sigma_{WX}$, thus relating the surrogate predictor to the true predictor. The auxiliary sample can be available under one of several scenarios in practice, which will be discussed in detail in Section 7. We will first (through Sections 3 to 6) focus on developments at the population level, and for this purpose it suffices to assume a matrix such as $\Sigma_{XW}$ is known, keeping in mind that it is to be obtained externally to the primary sample—either from the auxiliary data or from prior information.

Because $X$ is not observed, we use $\Sigma_{XW}$ to adjust the surrogate predictor $W$ to make it stochastically as close to $X$ as possible. As will soon be clear we can assume $E(X) = E(W) = 0$. In this case we adjust $W$ to $U = \Sigma_{XW}\Sigma_W^{-1}W$ (see, e.g., Carroll and Li [3]). Note that if $W$ is multivariate normal this is just the conditional expectation $E(X|W)$. Thus $U$ is the measurable function of $W$ closest to $X$ in terms of $L_2$ distance. If $W$ is not multivariate normal, then $U$ can simply be interpreted as linear regression of $X$ on $W$.

**3. Invariance of surrogate dimension reduction.** Recall that if there is a $p$ by $d$ matrix $\beta$, with $d \leq p$, such that (2) holds, then we call the column space of $\beta$ a dimension reduction space. See Li [19, 20] and Cook [6, 7]. Under very mild conditions, such as given in Cook [7], Section 6, the intersection of all dimension reduction spaces is again a dimension reduction space, which is then called the central space and is written as $\mathcal{S}_{Y|X}$. We will denote the dimension of $\mathcal{S}_{Y|X}$ by $q$. Note that $q \leq d$ for any $\beta$ satisfying (2). Similarly, we will denote the central space of $Y$ versus $U$ as $\mathcal{S}_{Y|U}$ and call it the surrogate central space. Our interests lie, of course, in the estimation of $\mathcal{S}_{Y|X}$, but $\mathcal{S}_{Y|U}$ is all that we can infer from the data. In this section we will establish the invariance law

$$(5) \qquad \mathcal{S}_{Y|U} = \mathcal{S}_{Y|X}$$

in the situation where $\Gamma$ is a $p$ by $p$ nonsingular matrix and $X$ and $\delta$ are multivariate normal.

We can assume without loss of generality that $E(X) = 0$ and $E(U) = 0$ because, for any $p$-dimensional vector $a$, $\mathcal{S}_{Y|X} = \mathcal{S}_{Y|(X-a)}$ and $\mathcal{S}_{Y|U} = \mathcal{S}_{Y|(U-a)}$. Since we will always assume $E(\delta) = 0$, $E(X) = E(U) = 0$ implies that $\gamma = 0$, and the measurement error model (1) reduces to

$$(6) \qquad W = \Gamma^T X + \delta.$$

The next lemma (and its variation) is the key to the whole development in this paper. It is also a fundamental fact about multivariate normal distributions that has been previously unknown. It will be applied to both exact and asymptotic distributions.

LEMMA 3.1. *Let $U_1^*, V_1^*$ be $r$-dimensional and $U_2^*, V_2^*$ be $s$-dimensional random vectors with $r + s = p$. Let*

$$V^* = \begin{pmatrix} V_1^* \\ V_2^* \end{pmatrix} \quad \text{and} \quad U^* = \begin{pmatrix} U_1^* \\ U_2^* \end{pmatrix},$$

*and let $Y$ be a random variable. Suppose:*

1. *$U_1^*, U_2^*, V_1^*, V_2^*$ are multivariate normal.*
2. *$U^* - V^* \perp\!\!\!\perp (V^*, Y)$.*

*Then:*

1. *If there is an $r$-dimensional multivariate normal random vector $V_3^*$ such that $V_3^* \perp\!\!\!\perp (V_1^* - V_3^*, V_2^*)$, and if $U_1^* \perp\!\!\!\perp U_2^*$, then $Y \perp\!\!\!\perp V^* | V_3^*$ implies $Y \perp\!\!\!\perp U^* | U_1^*$.*
2. *If there is an $r$-dimensional multivariate normal random vector $U_3^*$ such that $U_3^* \perp\!\!\!\perp (U_1^* - U_3^*, U_2^*)$, and if $V_1^* \perp\!\!\!\perp V_2^*$, then $Y \perp\!\!\!\perp U^* | U_3^*$ implies $Y \perp\!\!\!\perp V^* | V_1^*$.*

We should emphasize that despite its appearance the lemma is not symmetric for $U^*$ and $V^*$ because of assumption 2; note that we do not assume $V^* - U^* \perp\!\!\!\perp (U^*, Y)$. This is why the second assertion, though similar to the first, is not redundant. This asymmetry is intrinsic to the measurement error problem, where $U$ is a diffusion of $X$ but not conversely.

PROOF OF LEMMA 3.1. Write $U^*$ as $V^* + (U^* - V^*)$, and we have

(7) $$E(e^{it^T U^*} | Y) = E(e^{it^T V^*} e^{it^T (U^* - V^*)} | Y).$$

By assumption 2 we have $(U^* - V^*) \perp\!\!\!\perp V^* | Y$ and $U^* - V^* \perp\!\!\!\perp Y$. Hence the right-hand side reduces to

(8) $$E(e^{it^T V^*} | Y) E(e^{it^T (U^* - V^*)} | Y) = E(e^{it^T V^*} | Y) E(e^{it^T (U^* - V^*)}).$$

Assumption 2 also implies that $\Sigma_{U^*} = \text{var}(U^* - V^*) + \Sigma_{V^*}$, and hence that $U^* - V^* \sim N(0, \Sigma_{U^*} - \Sigma_{V^*})$. Thus the right-hand side of (8) further reduces to

$$E(e^{it^T V^*} | Y) e^{-(1/2) t^T (\Sigma_{U^*} - \Sigma_{V^*}) t}.$$

Substitute this into the right-hand side of (7) to obtain

(9) $$E(e^{it^T U^*} | Y) e^{(1/2) t^T \Sigma_{U^*} t} = E(e^{it^T V^*} | Y) e^{(1/2) t^T \Sigma_{V^*} t}.$$

Now suppose there is a $V_3^*$ such that $V_3^* \perp\!\!\!\perp (V_1^* - V_3^*, V_2^*)$ and $Y \perp\!\!\!\perp (V_1^*, V_2^*) | V_3^*$. The latter independence implies $Y \perp\!\!\!\perp (V_1^* - V_3^*, V_2^*) | V_3^*$ which, combined with the former independence, yields

$$(V_1^* - V_3^*, V_2^*) \perp\!\!\!\perp (V_3^*, Y) \Rightarrow (V_1^* - V_3^*, V_2^*) \perp\!\!\!\perp V_3^* | Y$$

and

$$(V_1^* - V_3^*, V_2^*) \perp\!\!\!\perp Y.$$

Hence

$$
\begin{aligned}
E(e^{it^T V^*}|Y) &= E(e^{it_1^T(V_1^* - V_3^*)} e^{it_2^T V_2^*} e^{it_1^T V_3^*}|Y) \\
&= E(e^{it_1^T(V_1^* - V_3^*)} e^{it_2^T V_2^*}) E(e^{it_1^T V_3^*}|Y).
\end{aligned}
\tag{10}
$$

Now let $W^* = ((V_1^* - V_3^*)^T, V_2^{*T})^T$. Because $V_3^* \perp\!\!\!\perp V_2^*$ we have

$$\Sigma_{W^*} = \begin{pmatrix} \Sigma_{V_1^* - V_3^*} & \Sigma_{V_1^* V_2^*} \\ \Sigma_{V_2^* V_1^*} & \Sigma_{V_2^*} \end{pmatrix} = \begin{pmatrix} \Sigma_{V_1^* - V_3^*} - \Sigma_{V_1^*} & 0 \\ 0 & 0 \end{pmatrix} + \Sigma_{V^*}. \tag{11}$$

In the meantime, because $W^*$ is multivariate normal we have

$$E(e^{it_1^T(V_1^* - V_3^*)} e^{it_2^T V_2^*}) = e^{-(1/2)t^T \Sigma_{W^*} t}. \tag{12}$$

Now combine (9) through (12) to obtain

$$E(e^{it^T U^*}|Y) e^{t^T \Sigma_{U^*} t} = E(e^{it_1^T V_3^*}|Y) e^{(1/2)t_1^T \Sigma_{V_1^*} t_1 - (1/2)t_1^T (\Sigma_{V_1^* - V_3^*}) t_1}.$$

Consequently the left-hand side does not depend on $t_2$; that is, we can take $t_2 = 0$ without changing it:

$$E(e^{it^T U^*}|Y) e^{t^T \Sigma_{U^*} t} = E(e^{it_1^T U_1^*}|Y) e^{t_1^T \Sigma_{U_1^*} t_1}.$$

This is equivalent to

$$E(e^{it^T U^*}|Y) = E(e^{it_1^T U_1^*}|Y) e^{-t^T \Sigma_{U^*} t + t_1^T \Sigma_{U_1^*} t_1} = E(e^{it_1^T U_1^*}|Y) e^{-t_2^T \Sigma_{U_2^*} t_2}, \tag{13}$$

where the second equality follows from the assumption $U_1^* \perp\!\!\!\perp U_2^*$. Multiply both sides by $e^{i\tau Y}$ and then take the expectation to obtain

$$E(e^{i\tau Y + it^T U^*}) = E[e^{i\tau Y} E(e^{it_1^T U_1^*}|Y)] e^{-t_2^T \Sigma_{U_2^*} t_2} = E(e^{i\tau Y} e^{it_1^T U_1^*}) e^{-t_2^T \Sigma_{U_2^*} t_2},$$

from which it follows that $(Y, U_1^*) \perp\!\!\!\perp U_2^*$, which implies that $Y \perp\!\!\!\perp U^*|U_1^*$.

The second assertion can be proved similarly. Following the same argument that leads to (10), we have

$$E(e^{it^T U^*}|Y) = E(e^{it_1^T U_3^*}|Y) e^{-(1/2)t_1^T (\Sigma_{U_1^* - U_3^*}) t_1 - (1/2)t_2^T \Sigma_{U_2^*} t_2}. \tag{14}$$

Now combine this relation with (9) and follow the proof of the first assertion to complete the proof. $\square$

We are now ready to establish the invariance relation (5).

THEOREM 3.1. *Suppose*:

1. $X \sim N(\mu_X, \Sigma_X)$, where $\Sigma_X > 0$.
2. $\delta \perp\!\!\!\perp (X, Y)$ and $\delta \sim N(0, \Sigma_\delta)$, where $\Sigma_\delta > 0$.

Then $\mathcal{S}_{Y|U} = \mathcal{S}_{Y|X}$.

PROOF. Assume, without loss of generality, that $E(X) = 0$ and $E(U) = 0$. We first show that $\mathcal{S}_{Y|U} \subseteq \mathcal{S}_{Y|X}$. Denote the dimension of $\mathcal{S}_{Y|X}$ by $q$, and let $\beta$ be a $p$ by $q$ matrix whose columns form a basis of $\mathcal{S}_{Y|X}$. Let $\zeta$ be a $p$ by $p-q$ matrix such that $\zeta^T \Sigma_U \beta = 0$ and such that the matrix $\eta = (\beta, \zeta)$ is full rank. Let

$$V_1^* = \beta^T \Sigma_U \Sigma_X^{-1} X, \qquad V_2^* = \zeta^T \Sigma_U \Sigma_X^{-1} X,$$
$$V_3^* = (\beta^T \Sigma_U \beta)(\beta^T \Sigma_X \beta)^{-1} \beta^T X,$$
$$U_1^* = \beta^T U, \quad U_2^* = \zeta^T U.$$

Then

$$\operatorname{cov}(U_1^*, U_2^*) = \beta^T \Sigma_U \zeta = 0,$$
$$\operatorname{cov}(V_3^*, V_2^*) = (\zeta^T \Sigma_U \beta)(\beta^T \Sigma_X \beta)^{-1}(\beta^T \Sigma_U \beta) = 0,$$
$$\operatorname{cov}(V_3^*, V_1^* - V_3^*) = (\beta^T \Sigma_U \beta)(\beta^T \Sigma_X \beta)^{-1}(\beta^T \Sigma_U \beta)$$
$$\qquad - (\beta^T \Sigma_U \beta)(\beta^T \Sigma_X \beta)^{-1}(\beta^T \Sigma_U \beta) = 0.$$

It follows that $U_1^* \perp\!\!\!\perp U_2^*$ and $V_3^* \perp\!\!\!\perp (V_1^* - V_3^*, V_2^*)$. In the meantime, by definition,

$$U^* - V^* = \eta^T(U - \Sigma_U \Sigma_X^{-1} X).$$

However, recall that

$$U = \Sigma_{XW} \Sigma_W^{-1} \Gamma^T X + \Sigma_{XW} \Sigma_W^{-1} \delta = \Sigma_{XW} \Sigma_W^{-1} \Sigma_{WX} \Sigma_X^{-1} X + \Sigma_{XW} \Sigma_W^{-1} \delta$$
(15)
$$= \Sigma_U \Sigma_X^{-1} X + \Sigma_{XW} \Sigma_W^{-1} \delta,$$

where the second equality holds because $\Sigma_{WX} = \Gamma^T \Sigma_X$, which follows from the independence $X \perp\!\!\!\perp \delta$ and the definition of $W$; the third equality holds because $\Sigma_U = \Sigma_{XW} \Sigma_W^{-1} \Sigma_{WX}$. Hence $U^* - V^* = \eta^T \Sigma_{XW} \Sigma_W^{-1} \delta$, which is independent of $(V^*, Y)$. Finally, we note that $V_3^*$ is a one-to-one function of $\beta^T X$ and $V^*$ is a function of $X$. So $Y \perp\!\!\!\perp V^* | V_3^*$. Thus, by the first assertion of Lemma 3.1, we have $U^* \perp\!\!\!\perp Y | U_1^* \Rightarrow U \perp\!\!\!\perp Y | \beta^T U \Rightarrow \mathcal{S}_{Y|U} \subseteq \mathcal{S}_{Y|X}$.

To prove $\mathcal{S}_{Y|X} \subseteq \mathcal{S}_{Y|U}$, let $\beta$ be a matrix whose columns are a basis of $\mathcal{S}_{Y|U}$, and $\zeta$ be such that the columns of $(\beta, \zeta)$ are a basis of $\mathbb{R}^p$ and $\zeta^T \Sigma_X \beta = 0$. Let

$$U_1^* = \beta^T \Sigma_X \Sigma_U^{-1} U, \qquad U_2^* = \zeta^T \Sigma_X \Sigma_U^{-1} U,$$
$$U_3^* = (\beta^T \Sigma_X \beta)(\beta^T \Sigma_U \beta)^{-1} \beta^T U,$$
$$V_1^* = \beta^T X, \qquad V_2^* = \zeta^T X.$$

Now follow the proof of $\mathcal{S}_{Y|U} \subseteq \mathcal{S}_{Y|X}$, but this time apply the second assertion of Lemma 3.1, to complete the proof. $\square$

The assumptions made in Theorem 3.1 are roughly equivalent to those made in Lue [22], Theorem 1, though our dimension reduction assumption, $Y \perp\!\!\!\perp X | \beta^T X$, is weaker than the corresponding assumption therein, which is $Y = g(\beta^T X, \varepsilon)$ where $\varepsilon \perp\!\!\!\perp X$—it is easy to show that the latter implies the former. For example, if $Y = g(\beta^T X, \varepsilon)$ where $\varepsilon \perp\!\!\!\perp X | \beta^T X$, then $Y \perp\!\!\!\perp X | \beta^T X$ still holds but $\varepsilon$ is no longer independent of $X$. Except for this dimension reduction assumption, our assumptions are stronger than those made in Carroll and Li [3], Theorem 2.1. However, our conclusion is stronger than those in both of these papers, in that it reveals the intrinsic invariance relation between dimension reduction spaces, not limited to any specific dimension reduction methods.

In the next example, we will give a visual demonstration of the invariance law.

EXAMPLE 3.1. Let $p = r = 6$, $X \sim N(0, I_p)$, $\varepsilon \sim N(0, \sigma_\varepsilon^2)$, $\delta \sim N(0, \sigma_\delta^2 I_p)$ and $\delta \perp\!\!\!\perp (X, \varepsilon)$. Consider the measurement-error regression model

$$(16) \qquad Y = 0.4(\beta_1^T X)^2 + 3\sin(\beta_2^T X/4) + \sigma_\varepsilon \varepsilon, \qquad W = \Gamma^T X + \delta,$$

where $\beta_1 = (1, 1, 1, 0, 0, 0)^T$ and $\beta_2 = (1, 0, 0, 0, 1, 3)^T$, and $\Gamma$ is a $p \times p$ matrix with diagonal elements equal to 1 and off-diagonal elements equal to 0.5. We take $\sigma_\varepsilon = 0.2$, $\sigma_\delta = 1/6$, and generate $(X_i, Y_i, W_i)$, $i = 1, \ldots, 400$, from this model.

In Figure 1 the left panels are the scatter plots of $Y$ versus $\beta_1^T X$ (upper) and $\beta_2^T X$ (lower) from a sample of 400 observed $(X, Y)$. The 3D shape of $Y$ versus $\beta^T X$ is roughly a $U$-shaped surface tilted upward in an orthogonal direction. The right panels are the scatter plots for $Y$ versus $\beta_1^T U$ and $\beta_2^T U$. As predicted by the invariance law, the directions $\beta_1$ and $\beta_2$, which are in $\mathcal{S}_{Y|X}$, also capture most of the variation of $Y$ in its relation to $U$, although the variation for the surrogate problem is larger than that for the original problem.

**4. The invariance law for arbitrary $\Gamma$.** We now turn to the general case where $\Gamma$ is a $p \times r$ matrix. We will assume $r \leq p$, which makes sense because otherwise there will be redundancy in the surrogate predictor $W$. In this case $W$ is of dimension $r$, but the adjusted surrogate predictor $U$ still has dimension $p$, with a singular variance matrix $\Sigma_U$ if $r < p$. We will assume that the column space of $\Gamma$ contains the dimension reduction space for $Y|X$ (which always holds if $\Gamma$ is a nonsingular square matrix). This is a very natural assumption—it means that we can have measurement error, but this error cannot be so erroneous as to erase part of the true regression parameter.

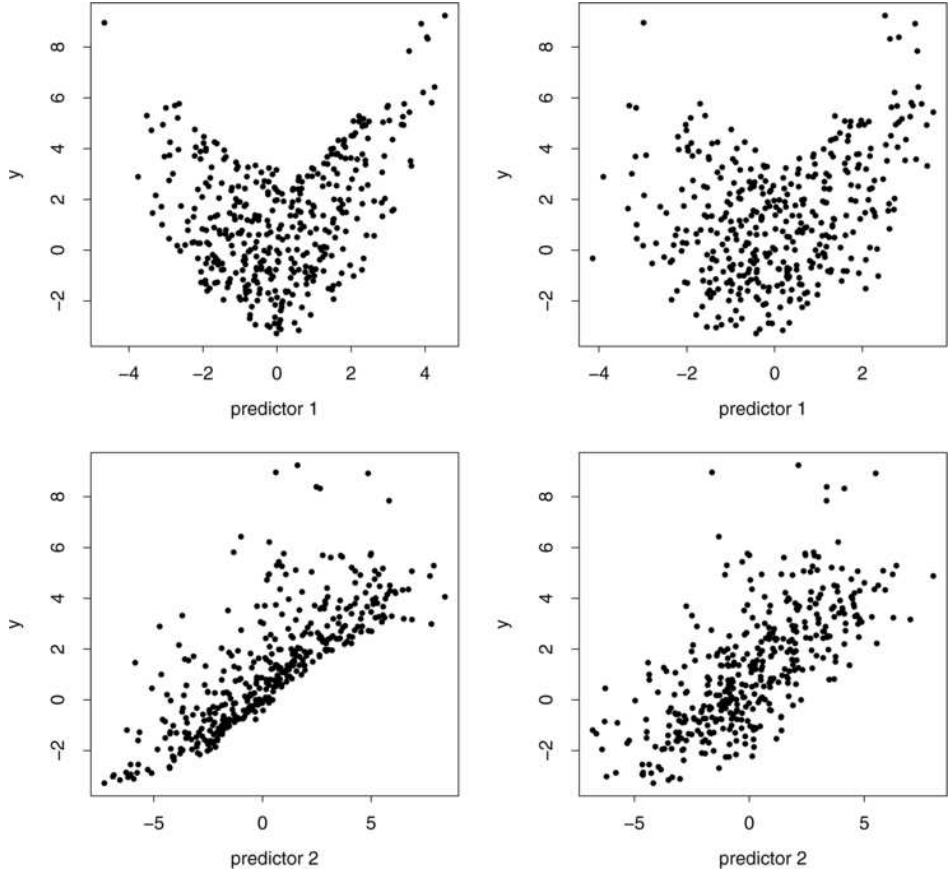

Fig. 1. *Original and surrogate dimension reduction spaces. Left panels: $Y$ versus $\beta_1^T X$ (upper) and $Y$ versus $\beta_2^T X$ (lower). Right panels: $Y$ versus $\beta_1^T U$ (upper) and $Y$ versus $\beta_2^T U$ (lower).*

THEOREM 4.1. *Suppose that $\Gamma$ in (1) is a $p$ by $r$ matrix with $r \leq p$, and that $\Gamma$ has rank $r$. Suppose that $\delta \sim N(0, \Sigma_\delta)$ with $\Sigma_\delta > 0$, $X \sim N(\mu_X, \Sigma_X)$ with $\Sigma_X > 0$, and $\delta \perp\!\!\!\perp (X, Y)$. Furthermore, suppose that $\mathcal{S}_{Y|X} \subseteq \mathrm{span}(\Gamma)$. Then $\mathcal{S}_{Y|U} = \mathcal{S}_{Y|X}$.*

PROOF. First we note that

$$\text{(17)} \qquad Y \perp\!\!\!\perp X | \Gamma^T X \quad \text{and} \quad Y \perp\!\!\!\perp U | \Gamma^T U.$$

The first relation follows directly from the assumption $\mathrm{span}(\beta) \subseteq \mathrm{span}(\Gamma)$. To prove the second relation, let $P_\Gamma(\Sigma_X) = \Gamma(\Gamma^T \Sigma_X \Gamma)^{-1} \Gamma^T \Sigma_X$ be the projection onto $\mathrm{span}(\Gamma)$ in terms of the inner product $\langle a, b \rangle = a^T \Sigma_X b$. Then

$$\mathrm{var}[(I - P_\Gamma(\Sigma_X))^T U] = [I - P_\Gamma(\Sigma_X)]^T \Sigma_U [I - P_\Gamma(\Sigma_X)].$$

However, we note that

$$[I - P_\Gamma(\Sigma_X)]^T \Sigma_U = (I - \Sigma_X \Gamma (\Gamma^T \Sigma_X \Gamma)^{-1} \Gamma^T) \Sigma_X \Gamma \Sigma_W^{-1} \Gamma^T \Sigma_X = 0.$$

Thus $\mathrm{var}[(I - P_\Gamma(\Sigma_X))^T U] = 0$, which implies $U = P_\Gamma^T(\Sigma_X) U$. That is, $U$ and $\Gamma^T U$ in fact generate the same $\sigma$-field, and hence the second relation in (17) must hold.

Next, by the definition of $U$ we have

$$U = \Sigma_X \Gamma \Sigma_W^{-1} (\Gamma^T X + \delta).$$

Multiply both sides of this equation from the left by $\Gamma^T$, to obtain

$$\Gamma^T U = \Gamma^T \Sigma_X \Gamma \Sigma_W^{-1} (\Gamma^T X + \delta).$$

Let $\widetilde{U} = \Gamma^T U$ and $\widetilde{X} = \Gamma^T X$. Then $\Sigma_{\widetilde{X}W} = \Gamma^T \Sigma_X \Gamma$ and $\Sigma_{\widetilde{U}} = \Gamma^T \Sigma_U \Gamma$. In this new coordinate system the above equation can be rewritten as

$$\widetilde{U} = \Sigma_{\widetilde{X}W} \Sigma_W^{-1} (\widetilde{X} + \delta).$$

Because (i) $\widetilde{X}$ has a multivariate normal distribution with $\Sigma_{\widetilde{X}} = \Gamma^T \Sigma_X \Gamma > 0$ and (ii) $\delta \perp\!\!\!\perp (\widetilde{X}, Y)$ and $\delta \sim N(0, \Sigma_\delta)$ with $\Sigma_\delta > 0$, we have, by Theorem 3.1, $\mathcal{S}_{Y|\widetilde{U}} = \mathcal{S}_{Y|\widetilde{X}}$.

Now let $q$ be the dimension of $\mathcal{S}_{Y|X}$ and suppose that $\beta$ is a $p$ by $q$ matrix whose columns form a basis of $\mathcal{S}_{Y|X}$. We note that $q \leq r$. Because $\mathrm{span}(\beta) \subseteq \mathrm{span}(\Gamma)$, there is an $r$ by $q$ matrix $\eta$ of rank $q$ such that $\beta = \Gamma \eta$. The following string of implications is evident:

$$Y \perp\!\!\!\perp X | \beta^T X \Rightarrow Y \perp\!\!\!\perp X | \eta^T \Gamma^T X \Rightarrow Y \perp\!\!\!\perp X | \eta^T \widetilde{X}$$
$$\Rightarrow Y \perp\!\!\!\perp \Gamma^T X | \eta^T \widetilde{X} \Rightarrow Y \perp\!\!\!\perp \widetilde{X} | \eta^T \widetilde{X}.$$

This means $\mathrm{span}(\eta)$ is a dimension reduction space for the problem $Y|\widetilde{X}$, and hence, because $\mathcal{S}_{Y|\widetilde{X}} = \mathcal{S}_{Y|\widetilde{U}}$, it must also be a dimension reduction space for the problem $Y|\widetilde{U}$. It follows that $Y \perp\!\!\!\perp \widetilde{U} | \eta^T \widetilde{U}$ or equivalently

$$(18) \qquad Y \perp\!\!\!\perp \Gamma^T U | \eta^T \Gamma^T U.$$

In the meantime, because $\Gamma^T U$ and $(\Gamma^T U, \eta^T \Gamma^T U)$ generate the same $\sigma$-field, the second relation in (17) implies

$$(19) \qquad Y \perp\!\!\!\perp U | (\Gamma^T U, \eta^T \Gamma^T U).$$

By Proposition 4.6 of Cook [7], relations (18) and (19) combined imply that

$$Y \perp\!\!\!\perp (U, \Gamma^T U) | \eta^T \Gamma^T U \Rightarrow Y \perp\!\!\!\perp U | \eta^T \Gamma^T U \Rightarrow Y \perp\!\!\!\perp U | \beta^T U,$$

from which it follows that $\mathcal{S}_{Y|U} \subseteq \mathcal{S}_{Y|X}$.

To show the reverse inclusion $\mathcal{S}_{Y|X} \subseteq \mathcal{S}_{Y|U}$, let $s$ be the dimension of $\mathcal{S}_{Y|U}$ and $\xi$ be a $p$ by $s$ matrix whose columns form a basis of $\mathcal{S}_{Y|U}$. By the second conditional independence in (17) we have $\text{span}(\xi) \subseteq \text{span}(\Gamma)$. Hence $s \le q$, and there is an $r$ by $s$ matrix $\zeta$ of rank $s$ such that $\xi = \Gamma \zeta$. Follow the proof of the first inclusion to show that

$$Y \perp\!\!\!\perp \Gamma^T X | \zeta^T \Gamma^T X.$$

In the meantime, because $\Gamma^T X$ and $(\Gamma^T X, \zeta^T \Gamma^T X)$ generate the same $\sigma$-field, the first conditional independence in (17) implies that

$$Y \perp\!\!\!\perp X | (\Gamma^T X, \zeta^T \Gamma^T X).$$

Now follow the proof of the first inclusion. □

**5. Invariance of surrogate dimension reduction for conditional moments.**
We now establish the invariance law between the central mean (or moment) spaces of the surrogate and the original dimension reduction problem. As briefly discussed in the Introduction, if there is a $p$ by $d$ matrix $\beta$ with $d \le p$ such that for $k = 1, 2, \ldots$,

$$(20) \qquad E(Y^k | X) = E(Y^k | \beta^T X),$$

then we call the column space of $\beta$ a dimension reduction space for the $k$th conditional moment. Similar to the previous case, the intersection of all such spaces again satisfies (20) under mild conditions. We call this intersection the $k$th central moment space, and denote it by $\mathcal{S}_{E(Y^k|X)}$. Let $\mathcal{S}_{E(Y^k|U)}$ be the $k$th central moment space for $Y$ versus $U$. The goal of this section is to establish the invariance relation

$$(21) \qquad \mathcal{S}_{E(Y^k|X)} = \mathcal{S}_{E(Y^k|U)}.$$

The next lemma parallels Lemma 3.1. Its proof will be given in the Appendix.

LEMMA 5.1. *Let $U_1^*, V_1^*, U_2^*, V_2^*, Y$ be as defined in Lemma* 3.1 *and suppose assumptions* 1 *and* 2 *therein are satisfied. Let $h(Y)$ be an integrable function of $Y$. Then*:

1. *If there is an $r$-dimensional multivariate normal random vector $V_3^*$ such that $V_3^* \perp\!\!\!\perp (V_1^* - V_3^*, V_2^*)$, and if $U_1^* \perp\!\!\!\perp U_2^*$, then $E[h(Y)|V^*, V_3^*] = E[h(Y)|V_3^*]$ implies $E[h(Y)|U^*] = E[h(Y)|U_1^*]$.*
2. *If there is an $r$-dimensional multivariate normal random vector $U_3^*$ such that $U_3^* \perp\!\!\!\perp (U_1^* - U_3^*, U_2^*)$, and if $V_1^* \perp\!\!\!\perp V_2^*$, then $E[h(Y)|U^*, U_3^*] = E[h(Y)|U_3^*]$ implies $E[h(Y)|V^*] = E[h(Y)|V_1^*]$.*

The next theorem establishes the invariance law (21).

THEOREM 5.1. *Suppose $k$ is any positive integer and*:

1. $X \sim N(\mu_X, \Sigma_X)$, *where* $\Sigma_X > 0$,
2. $\delta \perp\!\!\!\perp (X, Y)$ *and* $\delta \sim N(0, \Sigma_\delta)$, *where* $\Sigma_\delta > 0$,
3. $E(|Y|^k) < \infty$.

*Then* $\mathcal{S}_{E(Y^k|U)} = \mathcal{S}_{E(Y^k|X)}$.

The proof is similar to that of Theorem 3.1; the only difference is now we use Lemma 5.1 instead of Lemma 3.1. Evidently, we can follow the same steps in Section 4 to show that the assertion of Theorem 5.1 holds for general $\Gamma$. We state this generalization as the following corollary. The proof is omitted.

COROLLARY 5.1. *Suppose that $\Gamma$ in* (1) *is a $p$ by $r$ matrix with $r \leq p$, and that $\Gamma$ has rank $r$. Suppose that $\delta$ and $X$ are multivariate normal with $\Sigma_X > 0$, $\Sigma_\delta > 0$, $E(\delta) = 0$, and $\delta \perp\!\!\!\perp (X, Y)$. Suppose that $E(|Y|^k) < \infty$. Furthermore, suppose that $\mathcal{S}_{E(Y^k|X)} \subseteq \text{span}(\Gamma)$. Then $\mathcal{S}_{E(Y^k|U)} = \mathcal{S}_{E(Y^k|X)}$.*

**6. Approximate invariance for non-Gaussian predictors.** In this section we establish an approximate invariance law for arbitrary predictors. This is based on the fundamental result that, when the dimension $p$ is reasonably large, low-dimensional projections of the predictor are approximately multivariate normal. See Diaconis and Freedman [12] and Hall and Li [17]. Although this is a limiting behavior as $p \to \infty$, from our experience the approximate normality manifests for surprisingly small $p$. For example, a 1-dimensional projection of a 10-dimensional uniform distribution is virtually indistinguishable from a normal distribution. Thus the multivariate normality holds approximately in wide application. Intuitively, if we combine the exact invariance under normality, as we developed in the previous sections, and the approximate normality when $p$ is large, then we will arrive at an approximate invariance law for large $p$. This section is devoted to establishing this intuition as a fact.

We rewrite quantities such as $X, U, \delta, \beta$ in the previous sections as $X_p, U_p, \delta_p, \beta_p$. Let $\mathbb{S}^p$ denote the unit sphere in $\mathbb{R}^p : \{x \in \mathbb{R}^p : \|x\| = 1\}$, and $\text{Unif}(\mathbb{S}^p)$ denote the uniform distribution on $\mathbb{S}^p$. The result of Diaconis and Freedman [12] states that, if $\beta_p \sim \text{Unif}(\mathbb{S}^p)$, then, under mild conditions the conditional distribution of $\beta_p^T X_p | \beta_p$ converges *weakly in probability* (w.i.p.) to normal as $p \to \infty$. That is, the sequence of conditional characteristic functions $E(e^{it\beta_p^T X_p} | \beta_p)$ converges (pointwise) in probability to a normal characteristic function. Intuitively, this means when $p$ is large, the distribution of $\beta_p^T X_p$ is nearly normal for most $\beta_p$'s. Here, the parameter $\beta_p$ is treated as random to facilitate the notion of "most of $\beta_p$." We will adopt this assumption. In Diaconis and Freedman's development the $X_p$ is treated as nonrandom, but in our case $(\delta_p, X_p, Y)$ is random. In this context it makes sense to assume $\beta_p \perp\!\!\!\perp (Y, X_p, \delta_p)$, which would have been the case if the data $(Y, X_p, \delta_p)$ were treated as fixed. With $\beta_p$ being random, the dimension reduction relation should be stated as $Y \perp\!\!\!\perp X_p | (\beta_p^T X_p, \beta_p)$.

Our goal is to establish that, if $Y \perp\!\!\!\perp X_p | (\beta_p^T X_p, \beta_p)$, then, in an approximate sense $Y \perp\!\!\!\perp U_p | (\beta_p^T U_p, \beta_p)$, and vice versa. To do so we need to define an approximate version of conditional independence. Recall that, in the classical context when $p$ is fixed and $\beta$ is nonrandom, $Y \perp\!\!\!\perp U | \beta^T X$ if and only if $Y \perp\!\!\!\perp t^T U | \beta^T X$ for all $t \in \mathbb{R}^p$, as can be easily shown using characteristic functions. The definition of approximate conditional independence is analogous to the second statement.

DEFINITION 6.1. Let $\beta_p$ be a $p \times d$ dimensional random matrix whose columns are i.i.d. Unif($\mathbb{S}^p$) and $\beta_p \perp\!\!\!\perp (U_p, Y)$. We say that $Y$ and $U_p$ are asymptotically conditionally independent given $(\beta_p^T U_p, \beta_p)$, in terms of weak convergence in probability if, for any random vector $\zeta_p$ satisfying $\zeta_p \sim \text{Unif}(\mathbb{S}^p)$ and $\zeta_p \perp\!\!\!\perp (\beta_p, U_p, Y)$, the sequence $(Y, \beta_p^T U_p, \zeta_p^T U_p) | (\beta_p, \zeta_p)$ converges w.i.p. to $(Y^*, U^*, V^*)$ in which $Y^* \perp\!\!\!\perp V^* | U^*$. If this holds we write

$$Y \perp\!\!\!\perp U_p | (\beta_p^T U_p, \beta_p) \qquad \text{w.i.p. as } p \to \infty.$$

The following lemma gives further results regarding w.i.p. convergence that will be used in the later development. Its proof will be given in the [Appendix](#).

LEMMA 6.1. *Let $\{R_p\}$, $\{S_p\}$, $\{T_p\}$ and $\{\beta_p\}$ be sequences of random vectors in which the first three have dimensions not dependent on $p$. Then:*

1. *Let $R^*$ be a random vector with the same dimension as $R_p$ and denote by $\phi_p(t; \beta_p)$ and $\omega(t)$ the characteristic functions $E(e^{it^T R_p} | \beta_p)$ and $E(e^{it^T R^*})$, respectively. Then $R_p | \beta_p \to R^*$ w.i.p. if and only if, for each $t \in \mathbb{R}^p$,*

(22) $$E\phi_p(t; \beta_p) \to \omega(t), \qquad E|\phi_p(t; \beta_p)|^2 \to |\omega(t)|^2.$$

2. *If $(R_p, S_p, T_p) | \beta_p \to (R^*, S^*, T^*)$ w.i.p. and $R_p \perp\!\!\!\perp S_p | (T_p, \beta_p)$ for all $p$, then $R^* \perp\!\!\!\perp S^* | T^*$.*

Expression (22) is used as a sufficient condition for w.i.p. convergence in Diaconis and Freedman [12]; here we use it as a sufficient and necessary condition. In the next lemma, $\|\cdot\|_F$ will denote the Frobenius norm. Let $M_p$ be the random matrix $(U_p, X_p, \Sigma_{U_p} \Sigma_{X_p}^{-1} X_p)$ and $\widetilde{M}_p$ be an independent copy of $M_p$.

LEMMA 6.2. *If $\|\Sigma_{X_p}\|_F^2 = o(p^2)$, $\|\Sigma_{U_p}\|_F^2 = o(p^2)$ and $\|\Sigma_{U_p}\Sigma_{X_p}^{-1}\Sigma_{U_p}\|^2 = o(p^2)$, then $p^{-1}M_p^T \widetilde{M}_p = o_P(1)$.*

This will be proved in the Appendix. The convergence $p^{-1}M_p^T \widetilde{M}_p = o_P(1)$ was used by Diaconis and Freedman [12], as one of the two main assumptions [assumption (1.2)] in their development, but here we push this assumption back onto the structure of the covariance matrices. (More precisely, they used a parallel version of the convergence because they treat the data as a nonrandom sequence.) Conditions such as $\|\Sigma_{X_p}\|_F^2 = o(p^2)$ are quite mild. To provide intuition, recall that, if $\Sigma_p$ is a $p \times p$ matrix, and $\lambda_1, \ldots, \lambda_p$ are the eigenvalues of $\Sigma_p$ and $\lambda_{\max} = \max(\lambda_1, \ldots, \lambda_p)$, then

$$\|\Sigma_p\|_F^2 = \sum_{i=1}^p \lambda_i^2 \le p\lambda_{\max}^2.$$

Hence the condition $\|\Sigma_p\|_F^2 = o(p^2)$ will be satisfied if $\lambda_{\max} = o(\sqrt{p})$.

To streamline the assumptions, we make the following definition.

DEFINITION 6.2. We will say that a sequence of $p \times p$ matrices $\{\Sigma_p : p = 1, 2, \ldots\}$ is regular if $p^{-1}\operatorname{tr}(\Sigma_p) \to \sigma^2$ and $\|\Sigma_p\|_F^2 = o(p^2)$.

We now state the main result of this section.

THEOREM 6.1. *Suppose that $\beta_p$ is a $p \times d$ random matrix whose columns are i.i.d. $\operatorname{Unif}(\mathbb{S}^p)$. Suppose, furthermore, that:*

1. $\{\Sigma_{X_p}\}$, $\{\Sigma_{U_p}\}$ *and* $\{\Sigma_{U_p}\Sigma_{X_p}^{-1}\Sigma_{U_p}\}$ *are regular sequences with* $\sigma_X^2$, $\sigma_U^2$ *and* $\sigma_V^2$ *being the limits of their traces divided by $p$ as $p \to \infty$.*
2. $\delta_p \perp\!\!\!\perp (X_p, Y)$ *and* $\beta_p \perp\!\!\!\perp (X_p, Y, \delta_p)$.
3. $p^{-1}M_p^T M_p = E(p^{-1}M_p^T M_p) + o_P(1)$.

*If $Y \perp\!\!\!\perp X_p | (\beta_p^T X_p, \beta_p)$ and the conditional distribution of $Y|(\beta_p^T X_p = c, \beta_p)$ converges w.i.p. for each $c$, then $Y \perp\!\!\!\perp U_p | \beta_p^T U_p, \beta_p$ w.i.p. as $p \to \infty$.*

*If $Y \perp\!\!\!\perp U_p | (\beta_p^T U_p, \beta_p)$ and the conditional distribution of $Y|(\beta_p^T U_p = c, \beta_p)$ converges w.i.p. for each $c$, then $Y \perp\!\!\!\perp X_p | (\beta_p^T X_p, \beta_p)$ w.i.p. as $p \to \infty$.*

The condition that "the conditional distribution of $Y|(\beta_p^T X_p = c, \beta_p)$ converges w.i.p. for each $c$" means that the regression relation stabilizes as $p \to \infty$. Assumption 3 is parallel to the other of the two main assumptions in Diaconis and Freedman [12], Assumption 1.1. This is also quite mild: for example, it can be shown that if $X_p$ and $\delta_p$ are uniformly distributed on a ball $\{x \in \mathbb{R}^p : \|x\| \le \rho\}$ and if the covariance matrices involved satisfy some mild conditions, then assumption 3 is satisfied. For further discussion of this

assumption see Diaconis and Freedman [12], Section 3—though it is given in the context of nonrandom data, parallel conclusions can be drawn in our context.

PROOF OF THEOREM 6.1. For simplicity, we will only consider the case where $d = 1$; the proof of the general case is analogous and will be omitted. In this case $\beta_p \sim \text{Unif}(\mathbb{S}^p)$. Let $\zeta_p \sim \text{Unif}(\mathbb{S}^p)$ and $\zeta_p \perp\!\!\!\perp (\beta_p, X_p, \delta_p, Y)$. Following Diaconis and Freedman [12], we can equivalently assume $\beta_p \sim N(0, I_p/p)$ and $\zeta_p \sim N(0, I_p/p)$, because these distributions converge to $\text{Unif}(\mathbb{S}^p)$ as $p \to \infty$, and thus induce the same asymptotic effect as $\text{Unif}(\mathbb{S}^p)$. To summarize, we equivalently assume

$$\beta_p \sim N(0, I_p/p),$$
$$\zeta_p \sim N(0, I_p/p), \qquad \beta_p \perp\!\!\!\perp \zeta_p, \qquad (\beta_p, \zeta_p) \perp\!\!\!\perp (X_p, \delta_p, Y).$$

To prove the first assertion of the theorem, let

$$U_{1,p} = \beta_p^T U_p, \qquad U_{2,p} = \zeta_p^T U_p,$$
$$V_{1,p} = \beta_p^T \Sigma_{U_p} \Sigma_{X_p}^{-1} X_p, \qquad V_{2,p} = \zeta_p^T \Sigma_{U_p} \Sigma_{X_p}^{-1} X_p, \qquad V_{3,p} = (\sigma_U^2/\sigma_X^2) \beta_p^T X_p.$$

Our goal is to show that, as $p \to \infty$,

$$(Y, U_{1,p}, U_{2,p}, V_{1,p}, V_{2,p}, V_{3,p}) | (\beta_p, \zeta_p) \to (Y, U_1^*, U_2^*, V_1^*, V_2^*, V_3^*) \qquad \text{w.i.p.},$$

where $Y \perp\!\!\!\perp U_2^* | U_1^*$.

Let $(\beta_p, \zeta_p) = \eta_p$ and $L_p = (U_{1,p}, U_{2,p}, V_{1,p}, V_{2,p}, V_{3,p})^T$. Then

$$L_p = A(I_2 \otimes M_p^T) \text{vec}(\eta_p), \qquad \text{where } A = \begin{pmatrix} 1 & 0 & 0 & 0 & 0 & 0 \\ 0 & 0 & 0 & 1 & 0 & 0 \\ 0 & 0 & 1 & 0 & 0 & 0 \\ 0 & 0 & 0 & 0 & 0 & 1 \\ 0 & \dfrac{\sigma_U^2}{\sigma_X^2} & 0 & 0 & 0 & 0 \end{pmatrix}.$$

We first derive the w.i.p. limit of $(I_2 \otimes M_p^T) \text{vec}(\eta_p) | \eta_p$. Its (conditional) characteristic function is

$$\phi_p(t; \eta_p) = E(e^{it^T (I_2 \otimes M_p)^T \text{vec}(\eta_p)} | \eta_p), \qquad \text{where } t \in \mathbb{R}^{2p}.$$

Because $\text{vec}(\eta_p) \sim N(0, I_{2p}/p)$ and $\eta_p \perp\!\!\!\perp M_p$, we have

$$E[\phi_p(t; \eta_p)] = E[E(e^{it^T (I_2 \otimes M_p)^T \text{vec}(\eta_p)} | M_p)] = E(e^{-(1/(2p)) \|t^T (I_2 \otimes M_p)^T\|^2}).$$

By assumption 3 and assumption 1,

$$(23) \quad p^{-1} M_p^T M_p = p^{-1} E(M_p^T M_p) + o_P(1) \xrightarrow{P} \begin{pmatrix} \sigma_U^2 & \sigma_U^2 & \sigma_V^2 \\ \sigma_U^2 & \sigma_X^2 & \sigma_U^2 \\ \sigma_V^2 & \sigma_U^2 & \sigma_V^2 \end{pmatrix} \equiv \Gamma.$$

Thus, by the continuous mapping theorem (Billingsley [1], page 29),

$$(24) \quad e^{-(1/(2p))\|t^T(I_2 \otimes M_p)^T\|^2} \xrightarrow{P} e^{-(1/2)t^T(I_2 \otimes \Gamma)t}.$$

Because the sequence $\{e^{-(1/(2p))\|t^T(I_2 \otimes M_p)^T\|^2}\}$ is bounded, we have

$$(25) \quad E\phi_p(t; \eta_p) = E[e^{-(1/(2p))\|t^T(I_2 \otimes M_p)^T\|^2}] \to e^{-(1/2)t^T(I_2 \otimes \Gamma)t}.$$

In the meantime,

$$E[|\phi_p(t;\eta_p)|^2] = E[E(e^{it^T(I_2 \otimes M_p)^T \operatorname{vec}(\eta_p)}|\eta_p)E(e^{-it^T(I_2 \otimes M_p)^T \operatorname{vec}(\eta_p)}|\eta_p)].$$

If we let $\widetilde{M}_p$ be a copy of $M_p$ such that $\widetilde{M}_p \perp\!\!\!\perp (M_p, \eta_p)$, then the right-hand side can be rewritten as

$$E[E(e^{it^T[(I_2 \otimes M_p)^T - (I_2 \otimes \widetilde{M}_p)^T]\operatorname{vec}(\eta_p)}|M_p, \widetilde{M}_p)]$$
$$= E(e^{-(1/(2p))\|t^T[I_2 \otimes (M_p - \widetilde{M}_p)]^T\|^2}).$$

By Lemma 6.2 and convergence (23),

$$(M_p - \widetilde{M}_p)^T(M_p - \widetilde{M}_p) = M_p^T M_p - M_p^T \widetilde{M}_p - \widetilde{M}_p^T M_p + \widetilde{M}_p^T \widetilde{M}_p \xrightarrow{P} 2\Gamma.$$

Again, by the continuous mapping theorem,

$$e^{-(1/(2p))\|t^T[I_2 \otimes (M_p - \widetilde{M}_p)]^T\|^2} \xrightarrow{P} e^{-\|t^T(I_2 \otimes \Gamma)^T\|^2}.$$

Because the sequence $\{e^{-(1/(2p))\|t^T[I_2 \otimes (M_p - \widetilde{M}_p)]^T\|^2}\}$ is bounded, we have

$$(26) \quad E|\phi_p(t,s;\eta_p)|^2 = E[e^{-(1/(2p))\|t^T[I_2 \otimes (M_p - \widetilde{M}_p)]^T\|^2}] \to e^{-t^T(I_2 \otimes \Gamma)t}.$$

By part 1 of Lemma 6.1, (25) and (26),

$$(I_2 \otimes M_p)\operatorname{vec}(\eta_p)|\eta_p \to N(0, I_2 \otimes \Gamma) \quad \text{w.i.p. as } p \to \infty.$$

Thus the w.i.p. limit of $L_p|\eta_p$ is $N(0, A(I_2 \otimes \Gamma)A^T)$. By calculation,

$$\operatorname{cov}(U_1^*, U_2^*) = 0, \quad \operatorname{cov}(V_2^*, V_3^*) = 0,$$

$$\operatorname{cov}(V_1^* - V_3^*, V_3^*) = \operatorname{cov}(V_1^*, V_3^*) - \operatorname{cov}(V_3^*, V_3^*) = \frac{\sigma_U^4}{\sigma_X^2} - \frac{\sigma_U^4}{\sigma_X^2} = 0.$$

Hence, by multivariate normality we have $U_1^* \perp\!\!\!\perp U_2^*$ and $V_3^* \perp\!\!\!\perp (V_1^* - V_3^*, V_2^*)$. Also, recall from (15) that $U_p - \Sigma_{U_p}\Sigma_{X_p}^{-1}X_p = \Sigma_{X_pW_p}\Sigma_{W_p}^{-1}\delta_p$ and, by assumption 2, $\delta_p \perp\!\!\!\perp (X_p, Y)|\eta_p$. Consequently,

$$\eta_p^T U_p - \eta_p^T \Sigma_U \Sigma_X^{-1} X_p \perp\!\!\!\perp (X_p, Y)|\eta_p.$$

By part 2 of Lemma 6.1, $U^* - V^* \perp\!\!\!\perp (V^*, Y)$. So condition 2 of Lemma 3.1 is satisfied.

Let $L^*$ denote the w.i.p. limit of $L_p$. We now show that $(Y, L_p^T)^T \to (Y^*, L^{*T})^T$ w.i.p. for some random variable $Y^*$ and

(27) $$Y^* \perp\!\!\!\perp (V_1^*, V_2^*) | V_3^*.$$

Since, given $\eta_p$, $L_p$ is a function of $X_p$ and $\delta_p$, we have

$$E(e^{i\tau Y}|L_p, \eta_p) = E[E(e^{i\tau Y}|X_p, \delta_p, \eta_p)|L_p, \eta_p].$$

Since $Y \perp\!\!\!\perp X_p | (\beta_p^T X_p, \beta_p)$ and $\zeta_p \perp\!\!\!\perp (X_p, \delta_p, Y, \beta_p)$, we have $Y \perp\!\!\!\perp X_p | (\beta_p^T X_p, \eta_p)$. Also, since $\delta_p \perp\!\!\!\perp (X_p, Y, \eta_p)$, we have $Y \perp\!\!\!\perp \delta_p | (X_p, \eta_p)$. Hence

$$E(e^{i\tau Y}|X_p, \delta_p, \eta_p) = E(e^{i\tau Y}|X_p, \eta_p) = E(e^{i\tau Y}|V_{3,p}, \eta_p).$$

Thus $E(e^{i\tau Y}|L_p, \eta_p) = E(e^{i\tau Y}|V_{3,p}, \beta_p)$, or equivalently

(28) $$Y \perp\!\!\!\perp (L_p, \eta_p) | (V_{3,p}, \beta_p).$$

It follows that the conditional distribution of $Y|L_p, \eta_p$ is the same as the conditional distribution of $Y|(\beta_p^T X_p, \beta_p)$. Let $\mu(\cdot | c)$ be the w.i.p. limit of $Y|(\beta_p^T X_p = c, \beta_p)$, and draw the random variable $Y^*$ from $\mu(\cdot | V_1^*)$. Then $(Y, L_p)|\eta_p \to (Y^*, L^*)$ w.i.p. In the meantime (28) implies that $Y \perp\!\!\!\perp L_p | (V_{3,p}, \eta_p)$. Hence, by part 2 of Lemma 6.1, $Y \perp\!\!\!\perp L^* | V_3^*$, which implies (27). Thus all the conditions for assertion 1 of Lemma 3.1 are satisfied, and the first conclusion of this theorem holds.

The proof of the second assertion is similar, but this time let

$$U_{1,p} = \beta_p^T \Sigma_{X_p} \Sigma_{U_p}^{-1} U_p, \qquad U_{2,p} = \zeta_p^T \Sigma_{X_p} \Sigma_{U_p}^{-1} U_p, \qquad U_{3,p} = (\sigma_X^2/\sigma_U^2)\beta_p^T U_p,$$

$$V_{1,p} = \beta_p^T X_p, \qquad V_{2,p} = \zeta_p^T X_p$$

and use the second part of Lemma 3.1. The details are omitted. □

We now use a simulated example to demonstrate the approximate invariance law.

EXAMPLE 6.1. Still use the model in Example 3.1, but this time, assuming the distribution of $X$ is nonnormal, specified by

$$X_p = 3 Z_p \Phi(\|Z_p\|)/\|Z_p\|,$$

where $\Phi$ is the c.d.f. of the standard normal distribution and $Z_p \sim N(0, I_p)$. Thus, conditioning on each line passing through the origin, $X_p$ is uniformly distributed. Note that this is different from a uniform distribution over a ball in $\mathbb{R}^p$, but it is sufficiently nonnormal to illustrate our point. Figure 2 presents the scatter plots of $Y$ versus $\beta_1^T X_p$, $\beta_2^T X_p$ and the scatter plots of $Y$ versus $\beta_1^T U_p$, $\beta_2^T U_p$. We see that, even for $p$ as small as 6, $\mathcal{S}_{Y|U_p}$ already very much resembles $\mathcal{S}_{Y|X_p}$ for nonnormal predictors. In fact, we can hardly see any significant difference from Figure 1, where the invariance law holds exactly.

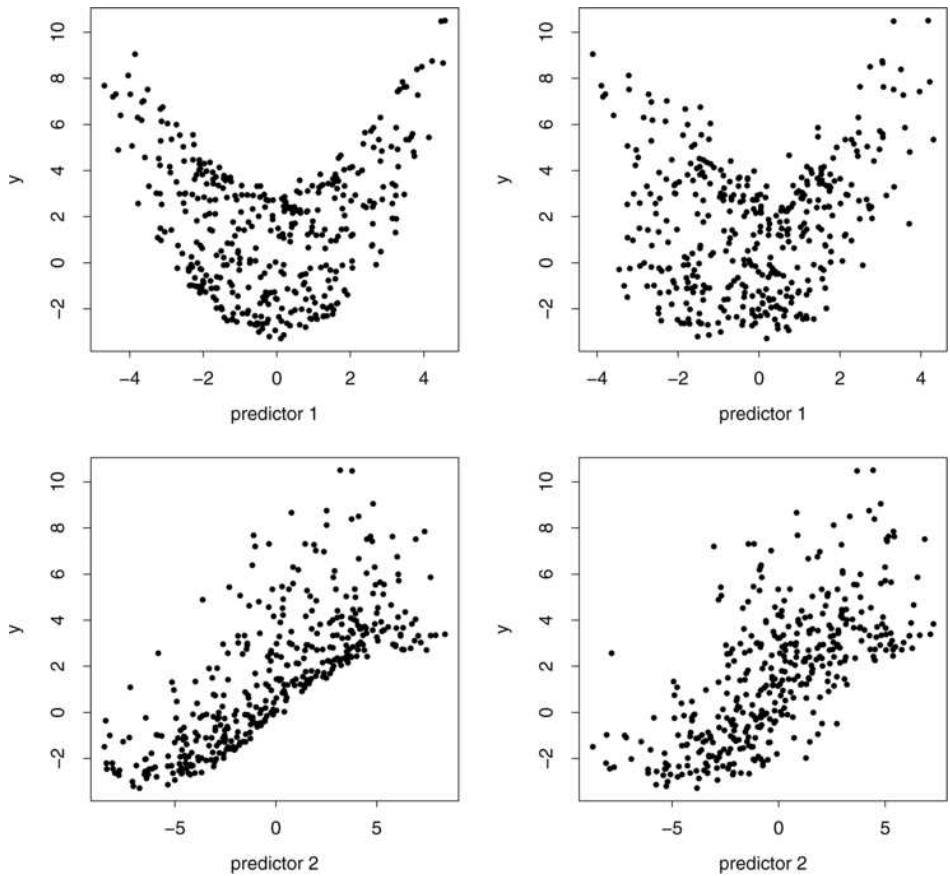

FIG. 2. *Comparison of the original and surrogate dimension reduction spaces for a non-normal predictor. Left panels: $Y$ versus $\beta_1^T X$ (upper) and $Y$ versus $\beta_2^T X$ (lower). Right panels: $Y$ versus $\beta_1^T U$ (upper) and $Y$ versus $\beta_2^T U$ (lower).*

Although we have only shown the asymptotic version of the invariance law (3) for nonsingular $p \times p$ dimensional $\Gamma$, similar extensions can be made for arbitrary $\Gamma$ (with $r \leq p$), as well as the invariance law (4). Because the development is completely analogous they will be omitted. Also notice that the assumptions for Theorem 6.1 impose no restriction on whether $X$ is a continuous random vector; thus the theorem also applies to discrete predictors—so long as its conditions are satisfied.

**7. Estimation procedure and convergence rate.** Having established the invariance laws at the population level, we now turn to the estimation procedure and the associated convergence rate for surrogate dimension reduction. Since we are no longer concerned with the limiting argument as $p \to \infty$, we will drop the subscript $p$ that indicates dimension. Instead the subscripts in $(X_i, Y_i)$ will now denote the $i$th case in the sample.

In the classical setting where measurement errors are absent, a dimension reduction estimator usually takes the following form. Let $(X_1, Y_1), \ldots, (X_n, Y_n)$ be an i.i.d. sample from $(X, Y)$. Let $F_{XY}$ be the distribution of $(X, Y)$, and $F_{n, XY}$ be the empirical distribution based on the sample. Let $\mathcal{F}$ be the class of all distributions of $(X, Y)$, and $\mathcal{G}$ be the class of all $p$ by $t$ matrices. Let $T : \mathcal{F} \to \mathcal{G}$ be a mapping from $\mathcal{F}$ to $\mathcal{G}$. Most of the existing dimension reduction methods, such as those described in the Introduction, have the form of such a mapping $T$, so chosen that (i) $\mathrm{span}(T(F_{XY})) \subseteq \mathcal{S}_{Y|X}$, and (ii) $T(F_{n,XY}) = T(F_{XY}) + \Delta_n$, where $\Delta_n$ is $o_p(1)$ or $O_p(1/\sqrt{n})$ depending on the estimators used. If these two conditions are satisfied with $\Delta_n = O_p(1/\sqrt{n})$ then we say that $T(F_{n,XY})$ is a $\sqrt{n}$-consistent estimator of $\mathcal{S}_{Y|X}$; if, in addition, (i) holds with equality then we say that $T(F_{n,XY})$ is a $\sqrt{n}$-exhaustive estimator of $\mathcal{S}_{Y|X}$. See Li, Zha and Chiaromonte [18].

The invariance law established in the previous sections tells us that we can apply a classical dimension reduction method to the adjusted surrogate predictor $U_1, \ldots, U_n$ and, if it satisfies (i) and (ii) for estimating $\mathcal{S}_{Y|U}$ (or $\mathcal{S}_{E(Y^k|U)}$), then it also satisfies these properties for estimating $\mathcal{S}_{Y|X}$ (or $\mathcal{S}_{E(Y^k|X)}$). Of course, here the adjusted surrogate predictor $U$ is not directly observed, because it contains such unknown parameters as $\Sigma_{XW}$ and $\Sigma_W$. However, these parameters can be estimated from an auxiliary sample that contains the information about the relation between $X$ and $W$. As discussed in Fuller [15] and Carroll and Li [3], in practice this auxiliary sample can be available under one of the several scenarios.

1. *Representative validation sample.* Under this scenario we observe, in addition to $(W_1, Y_1), \ldots, (W_n, Y_n)$, a validation sample

   $$(29) \qquad (X_{-1}, W_{-1}), \ldots, (X_{-m}, W_{-m}).$$

   We can use this auxiliary sample to estimate $\Sigma_{XW}$ and $\Sigma_W$ by the method of moments,

   $$\widehat{\Sigma}_{XW} = E_m[(X - \overline{X})(W - \overline{W})^T], \qquad \widehat{\Sigma}_W = E_m[(W - \overline{W})(W - \overline{W})^T],$$

   where $E_m$ denotes averaging over the representative validation sample (29). We can then use the estimates $\widehat{\Sigma}_{XW}$ and $\widehat{\Sigma}_W$ to adjust the surrogate predictor $W_i$ as $\widehat{U}_i = \widehat{\Sigma}_{XW} \widehat{\Sigma}_W^{-1} W_i$, and estimate $\mathcal{S}_{Y|X}$ by $T(F_{n,m,\widehat{U}Y})$. Here, $F_{n,m,\widehat{U}Y}$ is $F_{n,UY}$ with $U$ replaced by $\widehat{U}$; we have added the subscript $m$ to emphasize the dependence on $m$.

2. *Representative replication sample.* In this case we assume that $p = r$ and that $\Gamma$ is known, which, without loss of generality, can be taken as $I_p$. We have a separate sample where the true predictor $X_i$ is measured twice with measurement error. That is,

$$W_{ij} = \gamma + X_i + \delta_{ij}, \qquad j = 1, 2, \ i = -1, \ldots, -m, \tag{30}$$

where $\{\delta_{ij}\}$ are i.i.d. $N(0, \Sigma_\delta)$, $\{X_i\}$ are i.i.d. $N(0, \Sigma_X)$ and $\{\delta_{ij}\} \perp\!\!\!\perp \{(X_i, Y_i)\}$. From the replication sample we can extract information about $\Sigma_{XW}$ and $\Sigma_W$ as follows. It is easy to see that, for $i = -1, \ldots, -m$,

$$\mathrm{var}(W_{i1} - W_{i2}) = 2\Sigma_\delta, \qquad \mathrm{var}(W_{i1} + W_{i2}) = 4\Sigma_X + 2\Sigma_\delta,$$

from which we deduce

$$\Sigma_\delta = \mathrm{var}(W_{i1} - W_{i2})/2,$$
$$\Sigma_W = \mathrm{var}(W_{i1} + W_{i2})/4 + \mathrm{var}(W_{i1} - W_{i2})/4.$$

We can then estimate these two variance matrices by substituting in the right-hand side the moment estimates of $\mathrm{var}(W_{i1} + W_{i2})$ and $\mathrm{var}(W_{i1} - W_{i2})$ derived from the replication sample (30). Because in this case

$$\Sigma_{XW}\Sigma_W = I_p - \Sigma_\delta \Sigma_W^{-1},$$

we adjust the surrogate predictor $W_i$ as $\widehat{U}_i = (I_p - \widehat{\Sigma}_\delta \widehat{\Sigma}_W^{-1})W_i$.

A variation of the second sampling scheme appears in Fuller [15], which will be further discussed in Section 9 in conjunction with the analysis of a data set concerning managerial behavior.

Under the above schemes $\Sigma_{XW}$ and $\Sigma_W$ can be estimated by $\widehat{\Sigma}_{XW}$ and $\widehat{\Sigma}_W$ at the $\sqrt{m}$-rate, as can be easily derived from the central limit theorem. Hence $F_{n,m,\widehat{U}Y} = F_{n,UY} + O_p(1/\sqrt{m})$. Meanwhile, by the central limit theorem, we have $F_{n,UY} = F_{UY} + O_p(1/\sqrt{n})$. The dimension reduction estimator $T(F_{n,m,\widehat{U}Y})$ can be decomposed as

$$\begin{aligned}T(F_{n,m,\widehat{U}Y}) \\ = T(F_{UY}) + [T(F_{n,m,\widehat{U}Y}) - T(F_{n,UY})] + [T(F_{n,UY}) - T(F_{UY})].\end{aligned} \tag{31}$$

Usually the mapping $T$ is sufficiently smooth so that the second term is of order $O_p(1/\sqrt{m})$ and the third term is of the order $O_p(1/\sqrt{n})$. That is,

$$T(F_{n,m,\widehat{U}Y}) = T(F_{UY}) + O_p(1/\sqrt{m}) + O_p(1/\sqrt{n}). \tag{32}$$

While it is possible to impose a general smoothness condition on $T$ for the above relation in terms of Hadamard differentiability (Fernholz [14]), it is often simpler to verify (32) directly on an individual basis. The next example will illustrate how this can be done for a specific dimension reduction estimator.

EXAMPLE 7.1. Li, Zha and Chiaromonte [18] introduced contour regression (CR) which can be briefly described as follows. Let

$$H_X(c) = E[(\widetilde{X} - X)(\widetilde{X} - X)^T I(|\widetilde{Y} - Y| \le c)], \tag{33}$$

in which $(X,Y)$ and $(\widetilde{X}, \widetilde{Y})$ are independent random pairs distributed as $F_{XY}$, and $c > 0$ is a cutting point roughly corresponding to the width of a slice in SIR or SAVE. Let $v_{p-q+1}, \ldots, v_p$ be the eigenvectors of $\Sigma_X^{-1/2} H_X \Sigma_X^{-1/2}$ corresponding to its $q$ smallest eigenvalues. Then, under reasonably mild assumptions,

$$\mathrm{span}(\Sigma_X^{-1/2} v_p, \ldots, \Sigma_X^{-1/2} v_{p-q+1}) = \mathcal{S}_{Y|X}. \tag{34}$$

Thus, the mapping $T$ in this special case is defined by

$$T(F_{XY}) = (\Sigma_X^{-1/2} v_p, \ldots, \Sigma_X^{-1/2} v_{p-q+1}).$$

For estimation we replace $H_X$ and $\Sigma_X$ by their sample estimators,

$$\widehat{H}_X(c) = \binom{n}{2}^{-1} \sum_{(i,j)\in N} [(X_j - X_i)(X_j - X_i)^T I(|Y_j - Y_i| \le c)],$$
$$\widehat{\Sigma}_X = E_n[(X - \overline{X})(X - \overline{X})^T], \tag{35}$$

where, in the first equality, $N$ is the index set $\{(i,j): 1 \le j < i \le n\}$.

The motivation for introducing contour regression is to overcome the difficulties encountered by the classical methods. It is well known that if the relation $Y|X$ is symmetric about the origin then the population version of the SIR matrix is 0, and if $E(Y|X)$ is linear in $X$ then the population version of the pHd matrix is zero. In these cases, or in situations close to these cases, or in a combination of these cases, SIR and pHd tend not to yield accurate estimation of the dimension reduction space. Contour regression does not have this drawback because of the property (34).

In the context of measurement error regression the true predictor $X_i$ is to be replaced by $\widehat{U}_i$. For illustration, we adopt the first sampling scheme described above. Let $\widehat{\Sigma}_{W1}$ and $\widehat{\Sigma}_{W2}$ be the sample estimates of $\Sigma_W$ based on the primary sample $\{W_1, \ldots, W_n\}$ and the auxiliary sample $\{W_{-1}, \ldots, W_{-m}\}$, respectively. Let $\widehat{H}_{\widehat{U}}(c)$ be the matrix $\widehat{H}_X(c)$ in (35) with $X_i, X_j$ replaced by $\widehat{U}_i, \widehat{U}_j$. Then, it can be easily seen that

$$\widehat{H}_{\widehat{U}}(c) = \widehat{\Sigma}_{XW} \widehat{\Sigma}_{W2}^{-1} \widehat{H}_W(c) \widehat{\Sigma}_{W2}^{-1} \widehat{\Sigma}_{WX},$$
$$\widehat{\Sigma}_{\widehat{U}} = \widehat{\Sigma}_{XW} \widehat{\Sigma}_{W2}^{-1} \widehat{\Sigma}_{W1} \widehat{\Sigma}_{W2}^{-1} \widehat{\Sigma}_{WX}, \tag{36}$$

where in the first equality $\widehat{H}_W(c)$ is $\widehat{H}_X(c)$ in (35) with $X_i, X_j$ replaced by $W_i, W_j$. Because $\widehat{\Sigma}_{XW}$ and $\widehat{\Sigma}_{W2}$ are based on the auxiliary sample, they approximate $\Sigma_{XW}$ and $\Sigma_W$ at the $\sqrt{m}$-rate. Because $\widehat\Sigma_{W1}$ is based on the primary sample, it estimates $\Sigma_W$ at the $\sqrt{n}$-rate. Consequently, $\widehat\Sigma_{\widehat U} = \Sigma_U + O_p(1/\sqrt{m}) + O_p(1/\sqrt{n})$, which implies

$$\widehat\Sigma_{\widehat U}^{-1/2} = \Sigma_U^{-1/2} + O_p(1/\sqrt{n}) + O_p(1/\sqrt{m}).$$

In the meantime, it can be shown using the central limit theorem for $U$-statistics (see Li, Zha and Chiaromonte [18]) that $\widehat H_W(c) = H_W(c) + O_p(1/\sqrt{n})$, where $H_W(c)$ is $H_X(c)$ in (33) with $X, \widetilde X$ replaced by $W, \widetilde W$. Hence

$$\begin{aligned}(37)\quad \widehat H_{\widehat U}(c) &= \Sigma_{XW}\Sigma_W^{-1} H_W(c) \Sigma_W^{-1}\Sigma_{WX} + O_p(1/\sqrt{n}) + O_p(1/\sqrt{m}) \\ &= H_U(c) + O_p(1/\sqrt{n}) + O_p(1/\sqrt{m}).\end{aligned}$$

Combining (36) and (37) we see that

$$\widehat\Sigma_{\widehat U}^{-1/2} \widehat H_{\widehat U} \widehat\Sigma_{\widehat U}^{-1/2} = \Sigma_U^{-1/2} H_U \Sigma_U^{-1/2} + O_p(1/\sqrt{n}) + O_p(1/\sqrt{m}).$$

It follows that the eigenvectors $\hat v_{p-q+1}, \ldots, \hat v_p$ of the matrix on the left-hand side converge to the corresponding eigenvectors of the matrix on the right-hand side, $v_{p-q+1}, \ldots, v_p$, at the rate of $O_p(1/\sqrt{n}) + O_p(1/\sqrt{m})$, and consequently,

$$\widehat\Sigma_{\widehat U}^{-1/2}\hat v_i = \Sigma_U^{-1/2} v_i + O_p(1/\sqrt{n}) + O_p(1/\sqrt{m}).$$

Thus we have verified the convergence rate expressed in (32).

It is possible to use the general argument above to carry out asymptotic analysis for a surrogate dimension reduction estimator, in which both the rates according to $m$ and $n$ are involved. This can then be used to construct test statistics. But because of limited space this will be developed in future research. Special cases are available for SIR and OLS (Carroll and Li [3]) and for pHd (Lue [22]).

**8. Simulation studies.** As already mentioned, a practical impact of the invariance law is that it makes all dimension reduction methods accessible to the measurement error problem, thereby allowing us to tackle the difficult situations that the classical methods do not handle well. We now demonstrate this point by applying SIR, pHd and CR to simulated samples from the same model and comparing their performances.

We still use the model in Example 3.1, but this time take $\Gamma = I_p$. Under this model $\Sigma_W = \Sigma_X + \Sigma_\delta$ and $\Sigma_{XW} = \Sigma_X$. We take an auxiliary sample of size $m = 100$. The standard deviations $\sigma_\delta$ and $\sigma_\varepsilon$ are taken to be 0.2, 0.4, 0.6. For the auxiliary sample, we simulate $\{W_{ij} : j = 1, 2, i = -1, \ldots, -m\}$ according to the representative replication scheme described in Section 7.

TABLE 1
*Comparison of different surrogate dimension reduction methods*

| $\sigma_\varepsilon$ | $\sigma_\delta$ | **SIR** | **pHd** | **CR** |
|---|---|---|---|---|
|  | 0.2 | $1.26 \pm 0.63$ | $1.35 \pm 0.56$ | $0.12 \pm 0.07$ |
| 0.2 | 0.4 | $1.33 \pm 0.58$ | $1.56 \pm 0.51$ | $0.16 \pm 0.09$ |
|  | 0.6 | $1.46 \pm 0.53$ | $1.57 \pm 0.46$ | $0.32 \pm 0.21$ |
|  | 0.2 | $1.44 \pm 0.57$ | $1.41 \pm 0.51$ | $0.12 \pm 0.08$ |
| 0.4 | 0.4 | $1.34 \pm 0.57$ | $1.5 \pm 0.5$ | $0.20 \pm 0.13$ |
|  | 0.6 | $1.36 \pm 0.59$ | $1.62 \pm 0.48$ | $0.34 \pm 0.19$ |
|  | 0.2 | $1.36 \pm 0.62$ | $1.50 \pm 0.55$ | $0.13 \pm 0.08$ |
| 0.6 | 0.4 | $1.44 \pm 0.53$ | $1.53 \pm 0.49$ | $0.21 \pm 0.19$ |
|  | 0.6 | $1.48 \pm 0.53$ | $1.70 \pm 0.45$ | $0.32 \pm 0.18$ |

For each of the nine combinations of the values of $\sigma_\varepsilon$ and $\sigma_\delta$, 100 samples of $\{(X_i, Y_i)\}$ and $\{W_{ij}\}$ are generated according to the above specifications.

The estimation error of a dimension reduction method is measured by the following distance between two subspaces of $\mathbb{R}^p$. Let $\mathcal{S}_1$ and $\mathcal{S}_2$ be subspaces of $\mathbb{R}^p$, and $P_1$ and $P_2$ be the projections onto $\mathcal{S}_1$ and $\mathcal{S}_2$ with respect to the usual inner product $\langle a, b \rangle = a^T b$. We define the (squared) distance between $\mathcal{S}_1$ and $\mathcal{S}_2$ as

$$\rho(\mathcal{S}_1, \mathcal{S}_2) = \|P_1 - P_2\|^2,$$

where $\|\cdot\|$ is the Euclidean matrix norm. The same distance was used in Li, Zha and Chiaromonte [18], which is similar to the distance used in Xia et al. [25].

For SIR, the response is partitioned into eight slices, each having 50 observations. For CR, the cutting point $c$ is taken to be 0.5, which roughly amounts to including 12% of the $\binom{400}{2} = 79800$ vectors $\widehat{U}_i - \widehat{U}_j$ corresponding to the lowest increments in the response. The results are presented in Table 1. The symbol $a \pm b$ in the last three columns stands for mean and standard error of the distances $\rho(\widehat{\mathcal{S}}_{Y|X}, \mathcal{S}_{Y|X})$ over the 100 simulated samples. We can see that CR achieves significant improvement over SIR and pHd across all the combinations of $\sigma_\delta$ and $\sigma_\varepsilon$. This is because the regression model (16) contains a symmetric component in the $\beta_1$ direction, which SIR cannot estimate accurately, and a roughly monotone component in the $\beta_2$ direction, which pHd cannot estimate accurately. In contrast, CR accurately captures both component.

To provide further insights, we use one simulated sample to see the comparison among different estimators. Figure 3 compares the performance of SIR, pHd and CR in estimating $\mathcal{S}_{Y|U}$ (or $\mathcal{S}_{Y|X}$). We see that SIR gives a good estimate for $\beta_2$ but a poor estimate for $\beta_1$, the opposite of the case for pHd, but CR performs very well in estimating both components.

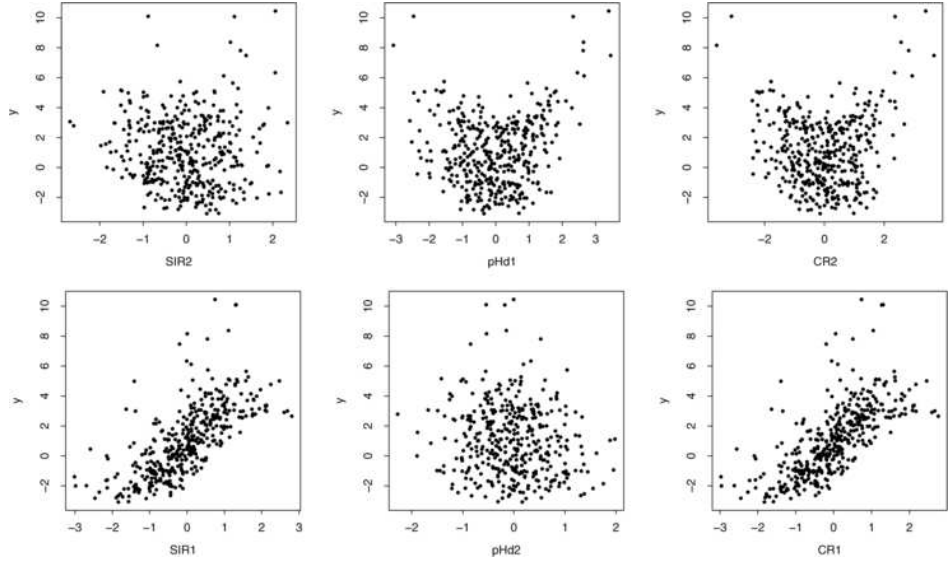

Fig. 3. *Surrogate dimension reduction by SIR, pHd and CR. Left panels: Y versus the second (upper) and the first (lower) predictors by SIR. Middle panels: Y versus first (upper) and second (lower) predictors by pHd. Right panels: the second (upper) and the first (lower) predictors by CR.*

**9. Analysis of a managerial role performance data.** In this section we apply surrogate dimension reduction to a role performance data set studied in Warren, White and Fuller [24] and Fuller [15]. To study managerial behavior, $n = 98$ managers of Iowa farmer cooperatives were randomly sampled. The response is the role performance of a manager. There are four predictors: $X_1$ (knowledge) is the knowledge of the economic phases of management directed toward profit-making, $X_2$ (value orientation) is the tendency to rationally evaluate means to an economic end and $X_3$ (role satisfaction) is gratification obtained

(training) is the amount of formal education. The predictors $X_1$, $X_2$ and $X_3$, and the response $Y$ are measured with questionnaires filled out by the managers, and contain measurement errors. The amount of formal education, $X_4$, is measured without error. Through dimension reduction of this data set we wish to see if the linear assumption in Fuller [15] is reasonable, if there are linear combinations of the predictors other than those obtained from the linear model that significantly affect the role performance, and how different surrogate dimension reduction methods perform and compare in a real data set.

The sampling scheme is a variation of the second one described in Section 7. A split halves design of the questionnaires yields two observations on each item with measurement error for each manager, say

$$(V_{i1}, V_{i2}), \qquad (W_{i11}, W_{i12}), \ldots, (W_{i31}, W_{i32}), \qquad i = 1, \ldots, n,$$

where $(V_{i1}, V_{i2})$ are the two measurements of $Y_i$ and $(W_{ij1}, W_{ij2})$, $j = 1, 2, 3$, are the two measurements of $X_{ij}$. The surrogate predictors for $X_{ij}$, $j = 1, 2, 3$, are taken to be the averages $W_{ij} = (W_{ij1} + W_{ij2})/2$. Similarly, the surrogate response for $Y_i$ is taken to be $V_i = (V_{i1} + V_{i2})/2$. Since the measurement error in $Y$ does not change the regression model, we can treat $V$ as the true response $Y$. As in Fuller [15], we assume:

1. $V_{ik} = Y_i + \xi_{ik}$, $W_{ijk} = X_{ij} + \eta_{ijk}$, $i = 1, \ldots, n$, $j = 1, 2, 3$, $k = 1, 2$, where
$$\{(\xi_{ik}, \eta_{i1k}, \eta_{i2k}, \eta_{i3k}) : i = 1, \ldots, n, k = 1, 2\}$$
$$\perp\!\!\!\perp \{(Y_i, X_{i1}, \ldots, X_{i4}) : i = 1, \ldots, n\}.$$

2. The random vectors $\{(\xi_{ik}, \eta_{i1k}, \eta_{i2k}, \eta_{i3k}) : i = 1, \ldots, n\}$ are i.i.d. 4-dimensional normal with mean 0 and variance matrix
$$\operatorname{diag}(\sigma_\xi^2, \sigma_{\eta,1}^2, \sigma_{\eta,2}^2, \sigma_{\eta,3}^2).$$

3. $\{(X_{i1}, \ldots, X_{i4}) : n = 1, \ldots, n\}$ are i.i.d. $N(\mu_X, \Sigma_X)$.

From these assumptions it is easy to see that, for $j = 1, \ldots, 4$ and $i = 1, \ldots, n$, $W_{ij} = X_{ij} + \delta_{ij}$, where

$$\{(\delta_{i1}, \ldots, \delta_{i4}) : i = 1, \ldots, n\} \perp\!\!\!\perp \{(X_{i1}, \ldots, X_{i4}, Y_i) : i = 1, \ldots, n\}.$$

and $\{(\delta_{i1}, \ldots, \delta_{i4}) : i = 1, \ldots, n\}$ are i.i.d. mean 0 and variance matrix

$$\Sigma_\delta = \operatorname{diag}(\sigma_{\delta,1}^2, \sigma_{\delta,2}^2, \sigma_{\delta,3}^2, 0), \qquad \sigma_{\delta,j}^2 = \tfrac{1}{4} \operatorname{var}(W_{ij1} - W_{ij2}),$$
(38)
$$j = 1, 2, 3.$$

Thus, at the population level, our measurement error model can be summarized as

$$W = X + \delta, \qquad \delta \perp\!\!\!\perp (X, V), \qquad \delta \sim N(0, \Sigma_\delta), \qquad X \sim N(\mu_X, \Sigma_X),$$

where $\Sigma_\delta$ is given by (38). Note that, unlike in Fuller [15], no assumption is imposed on the relation between $Y$ and $X$.

The measurement error variance $\Sigma_\delta$ is estimated using the moment estimator of (38) based on the sample $\{(W_{ij1}, W_{ij2}) : i = 1, \ldots, 98, j = 1, \ldots, 3\}$, as

$$\widehat{\Sigma}_\delta = \operatorname{diag}(0.0203, 0.0438, 0.0180, 0).$$

The variance matrix $\Sigma_W$ of $W_i = (W_{i1}, \ldots, W_{i4})^T$ is estimated from the sample $\{W_i : i = 1, \ldots, 98\}$ as

$$\begin{pmatrix} 0.0520 & 0.0280 & 0.0044 & 0.0192 \\ 0.0280 & 0.1212 & -0.0063 & 0.0353 \\ 0.0044 & -0.0063 & 0.0901 & -0.0066 \\ 0.0192 & 0.0353 & -0.0066 & 0.0946 \end{pmatrix}.$$

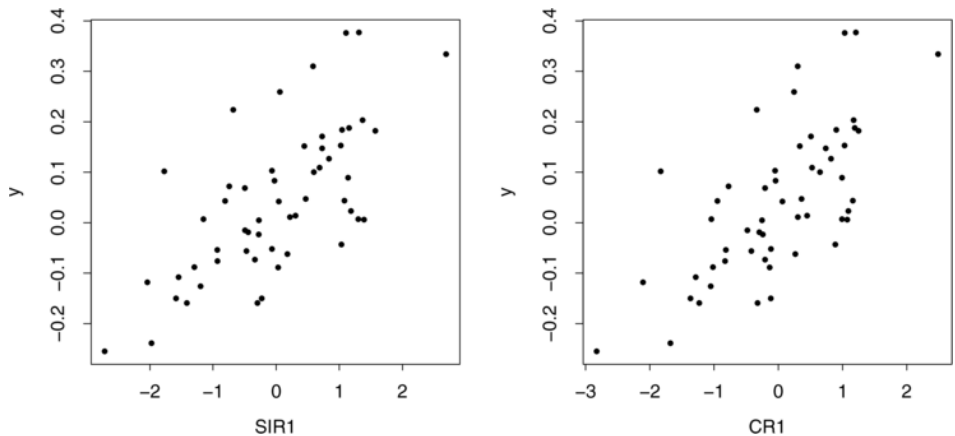

FIG. 4. *Scatter plots of role performance versus the first predictor by SIR (left) and versus the first predictor by CR (right).*

The correction coefficients $\Sigma_{XW}\Sigma_W^{-1}$ are then estimated by $I_4 - \widehat{\Sigma}_\delta \widehat{\Sigma}_W^{-1}$, and the surrogate predictor $W_i$ corrected as $\widehat{U}_i = (I_4 - \widehat{\Sigma}_\delta \widehat{\Sigma}_W^{-1})W_i$.

We apply SIR and contour regression to the surrogate regression problem of $V_i$ versus $\widehat{U}_i$. As in Fuller [15], a portion of the data (55 out of 98 subjects) will be presented. For SIR, the responses of 55 subjects are divided into five slices, each having 11 subjects. For CR, the cutting point $c$ is taken to be 0.1, which amounts to including 552 of the $\binom{55}{2} = 1485$ (roughly 37%) differences $\widehat{U}_i - \widehat{U}_j$ corresponding to the lowest increments $|Y_i - Y_j|$ of the response. In fact, varying the number of slices for SIR or the cutting point $c$ for CR within a reasonable range does not seem to have a serious effect on their performance.

Figure 4 shows the scatter plots of $Y$ versus the first predictors calculated by SIR (left) and CR (right). None of the scatter plots of $Y$ versus the second predictors by SIR and CR shows any discernible pattern, and so they are not presented. Because there is no $U$-shaped component in the data, the accuracy of SIR and CR are comparable. These plots show that the linear model postulated in Warren, White and Fuller [24] and Fuller [15] does fit this data, and there do not appear to be other linear combinations of the predictors that significantly affect the role performance. The directions obtained from CR, SIR, and that using the maximum likelihood estimator for a linear model, are presented in Table 2 (the vectors are rescaled to have lengths 1).

Note that SIR, as applied to the adjusted surrogate predictor $\widehat{U}$, is the estimator proposed in Carroll and Li [3]. We can see that for this data set the three methods are more or less consistent, though CR gives more weight to past education than the other methods. Of course, the significance of these parameters should be determined by a formal test based on the relevant asymptotic distributions. Such asymptotic results are available for pHd and SIR, and are under development for CR.

## APPENDIX

PROOF OF LEMMA 5.1. Note that the conditions for equality (9) are still satisfied. Multiply both sides of (9) by $h(Y)$ and then take expectation, to obtain

$$E[h(Y)e^{it^T U^*}]e^{(1/2)t^T \Sigma_{U^*} t} = E[h(Y)e^{it^T V^*}|Y]e^{(1/2)t^T \Sigma_{V^*} t}.$$

To prove the first assertion, suppose there is a $V_3^*$ such that $V_3^* \perp\!\!\!\perp (V_1^* - V_3^*, V_2^*)$ and $E[h(Y)|V^*, V_3^*] = E[h(Y)|V_3^*]$. Then

$$\begin{aligned}
E(h(Y)e^{it^T V^*}) &= E(E[h(Y)|V^*, V_3^*]e^{it_1^T(V_1^* - V_3^*)}e^{it_2^T V_2^*}e^{it_1^T V_3^*}) \\
&= E(E[h(Y)|V_3^*]e^{it_1^T(V_1^* - V_3^*)}e^{it_2^T V_2^*}e^{it_1^T V_3^*}) \\
&= E(e^{it_1^T(V_1^* - V_3^*)}e^{it_2^T V_2^*})E(h(Y)e^{it_1^T V_3^*}).
\end{aligned}$$

Follow the steps that lead to (13) in the proof of Lemma 3.1 to obtain

$$E(h(Y)e^{it^T U^*}) = E(h(Y)e^{it_1^T U_1^*})e^{-t_2^T \Sigma_{U_2^*} t_2}.$$

The equation can be rewritten as

$$E(E[h(Y)|U^*]e^{it^T U^*}) = E(E[h(Y)|U_1^*]e^{it_1^T U_1^*})e^{-t_2^T \Sigma_{U_2^*} t_2}.$$

Because $U_1^* \perp\!\!\!\perp U_2^*$, the right-hand side is

$$E(E[h(Y)|U_1^*]e^{it_1^T U_1^* + it_2^T U_2^*}) = E(E[h(Y)|U_1^*]e^{it^T U^*}).$$

It follows that

$$E(\{E[h(Y)|U^*] - E[h(Y)|U_1^*]\}e^{it^T U^*}) = 0$$

for all $t$. In other words, the Fourier transform of $E[h(Y)|U^*] - E[h(Y)|U_1^*]$ is zero. Thus $E[h(Y)|U^*] - E[h(Y)|U_1^*] = 0$ almost surely.

The proof of the second assertion will be omitted. □

TABLE 2

|  | $\hat{\beta}_1$ | $\hat{\beta}_2$ | $\hat{\beta}_3$ | $\hat{\beta}_4$ |
|---|---|---|---|---|
| Fuller [15] | 0.881 | 0.365 | 0.286 | 0.098 |
| SIR (Carroll and Li) | 0.952 | 0.219 | 0.187 | 0.102 |
| CR | 0.935 | 0.291 | 0.126 | 0.159 |

PROOF OF LEMMA 6.1. 1. That (22) implies $\phi_p(t;\beta_p) \xrightarrow{P} \omega(t)$ is well known. Now suppose $\phi_p(t;\beta_p) \xrightarrow{P} \omega(t)$. Then $|\phi_p(t;\beta_p)|^2 \xrightarrow{P} |\omega(t)|^2$. Because both $\phi_p(t;\beta_p)$ and $|\phi_p(t;\beta_p)|^2$ are bounded, (22) holds.

2. Because $R_p \perp\!\!\!\perp S_p | T_p, \beta_p$, we have

$$\begin{aligned}
&E(e^{it^T R^*} e^{iu^T S^*} | T^*) \\
&= E(e^{it^T R^*} | T^*) E(e^{iu^T S^*} | T^*) \\
&\quad + [E(e^{it^T R_p} | T_p, \beta_p) E(e^{iu^T S_p} | T_p, \beta_p) - E(e^{it^T R^*} | T^*) E(e^{iu^T S^*} | T^*)] \\
&\quad + [E(e^{it^T R^*} e^{iu^T S^*} | T^*) - E(e^{it^T R_p} e^{iu^T S_p} | T_p, \beta_p)].
\end{aligned}$$

Because $(R_p, S_p, T_p)|\beta_p \to (R^*, S^*, T^*)$ w.i.p., the last two terms on the right-hand side are $o_P(1)$. Hence the left-hand side equals the first term on the right-hand side because the former is a nonrandom quantity independent of $p$. □

PROOF OF LEMMA 6.2. It suffices to show that, if $A_p$ and $B_p$ are regular sequences of random vectors in $\mathbb{R}^p$ such that $A_p \perp\!\!\!\perp B_p$ and $E(A_p) = 0$, $E(B_p) = 0$, then $p^{-1} A_p^T B_p = o_P(1)$. If this is true then we can take $A_p = X_p, U_p$, or $\Sigma_U \Sigma_X^{-1} X_p$ and $B_p = \widetilde{X}_p, \widetilde{U}_p$, or $\Sigma_U \Sigma_X^{-1} \widetilde{X}_p$ to prove the desired equality.

By Chebyshev's inequality,

$$(39) \qquad P(p^{-1}|A_p^T B_p| > \varepsilon) \le \frac{1}{\varepsilon^2 p^2} E(A_p^T B_p)^2.$$

The expectation on the right-hand side is

$$E\left(\sum_{i=1}^p A_p^i B_p^i\right)^2 = \sum_{i=1}^p \sum_{j=1}^p E(A_p^i B_p^i A_p^j B_p^j)$$

$$= \sum_{i=1}^p \sum_{j=1}^p \Sigma_{A,ij} \Sigma_{B,ij} \le \|\Sigma_A\|_F \|\Sigma_B\|_F,$$

where the inequality is from the Cauchy–Schwarz inequality. By assumption, $\|\Sigma_A\|_F = o(p)$ and $\|\Sigma_B\| = o(p)$. So the right-hand side of (39) converges to 0 as $p \to \infty$. Hence $A_p^T B_p = o_P(1)$. □

**Acknowledgments.** We are grateful to two referees and an Associate Editor who, along with other insightful comments, suggested we consider the situations where the predictor and measurement error have non-Gaussian distributions, which led to the development in Section 6.

DEPARTMENT OF STATISTICS
326 THOMAS BUILDING
PENNSYLVANIA STATE UNIVERSITY
UNIVERSITY PARK, PENNSYLVANIA 16802
USA
E-MAIL: bing@stat.psu.edu

DEPARTMENT OF STATISTICS
FRANKLIN COLLEGE
107 STATISTICS BUILDING
UNIVERSITY OF GEORGIA
ATHENS, GEORGIA 30602-1952
USA
E-MAIL: xryin@stat.uga.edu